# A scalable multilevel domain decomposition preconditioner with a subspace-based coarsening algorithm for the neutron transport calculations

Fande Kong[†*], Yaqi Wang[‡], Derek R. Gaston[†], Alexander D. Lindsay[†], Cody J. Permann[†], Richard C. Martineau[†]

[†]*Computational Frameworks, Idaho National Laboratory, P.O. Box 1625, Idaho Falls, ID 83415-3840, USA*

[‡]*Nuclear Engineering Methods Development, Idaho National Laboratory, P.O. Box 1625, Idaho Falls, ID 83415, USA*

SUMMARY

The multigroup neutron transport equations has been widely used to study the interactions of neutrons with their background materials in nuclear reactors. High-resolution simulations of the multigroup neutron transport equations using modern supercomputers require the development of scalable parallel solving techniques. In this paper, we study a scalable transport method for solving the algebraic system arising from the discretization of the multigroup neutron transport equations. The proposed transport method consists of a fully coupled Newton solver for the generalized eigenvalue problems and GMRES together with a novel multilevel domain decomposition preconditioner for the Jacobian system. The multilevel preconditioner has been successfully used for many problems, but the construction of coarse spaces for certain problems, especially for unstructured mesh problems, is expensive and often unscalable. We introduce a new subspace-based coarsening algorithm to address this issue by exploring the structure of the matrix in the discretized version of the neutron transport problems. We numerically demonstrate that the proposed transport method is highly scalable with more than 10,000 processor cores for the 3D C5G7 benchmark problem on unstructured meshes with billions of unknowns. Compared with the traditional multilevel domain decomposition method, the new approach equipped with the subspace-based coarsening algorithm is much faster on the construction of coarse spaces.

## 1. INTRODUCTION

The multigroup neutron transport equations studies the interaction of neutrons with their background materials, and it has been widely used in nuclear reactor simulations. Finding an approximate numerical solution to the multigroup neutron transport equations remains as one of the most computationally expensive problems since the equations is defined on a high-dimensional space phase (3D spatial space, 2D angle, 1D energy) [1, 2]. The system of algebraic equations arising from the discretization of the multigroup neutron transport equations for completely describing the physical systems of interest is often very large,

*Correspondence to: (fande.kong@inl.gov) Computational Frameworks, Idaho National Laboratory, P.O. Box 1625, Idaho Falls, ID 83415, USA

and thus a scalable parallel solver that is capable of using thousands of processor cores is required. However, designing and developing a scalable parallel solver for the neutron transport equations is still a challenging task when the number of processor cores is large since there often exist a few algorithm components that are inherently sequential and these algorithm components will dominate the overall calculations when more and more processor cores are used. In this paper, we propose a scalable parallel multilevel domain decomposition preconditioner with a subspace-based coarsening algorithm to address this issue by exploring the structure of the matrix in the discretized version of the multigroup transport equations.

To place this work in the context of previous research, we present a brief literature review. A more complete literature review of the solving techniques used in the neutron transport simulations can be found in [3]. The development of efficient transport methods has been an active area of research for a couple of decades. The simplest and oldest method is known as transport *sweeps*, where a Gauss-Seidel iteration is implemented for energy groups and angular directions, and within an angular direction a small element matrix is inverted locally when sweeping through the mesh for computing the inversion of the collision and the streaming terms. In the transport sweeps, a source iteration (Richardson iteration) is often employed as the outer solver for the state-steady or transient simulations and an inverse power iteration is used as the outer loop for the eigenvalue calculations. For diffiusive problems, transport sweeps are slow to converge and require an acceleration technique. The *synthetic* method [4] and the *Diffusion Synthetic Acceleration (DSA)* method [5] are popular methods for such a purpose. In DSA, the diffusion equations is employed as a low-order operator to accelerate transport sweeps. The transport sweeps together with DSA requires inherently sequential operations that are not ideally suited to parallel computing [6, 7]. There are many works [8, 9] that try to overcome this difficulty. In fact, the combination of transport sweeps together with a low order diffusion-based acceleration can be viewed as a two-grid/two-level scheme in angle. Over the past decade much research has been done on applying multigrid/multilevel methods to the transport equations. In [10] for instance, the authors study angular multigrid methods that go beyond the two-grid/two-level scheme of DSA. In [7, 11] the spatial multigrid methods are employed for different problems. The multigrid/multilevel methods has achieved a great success for the neutron transport simulations since it shows a great convergence for many problems and is also well suited for parallel computing. But in the multigrid/multilevel methods, it is well-known that the construction of coarse spaces is often inherently sequential and not perfectly scalable. The setup time of the multigrid/multilevel methods may not be ideally decreased when the number of processor cores is increased.

To resolve this issue, we propose a highly scalable, fully coupled parallel transport method, where, instead of an inverse power iteration, a Jacobian-free Newton-based eigenvalue solver is employed for computing the fundamental mode of the neutron transport equations by reforming the generalized eigenvalue problem as a nonlinear system of equations. During each Newton iteration, the Jacobian system is calculated using GMRES [12] that requires an efficient preconditioner to speedup the convergence. We employ a multilevel domain decomposition method as a preconditioner with a novel approach for constructing coarse spaces. The new approach is much cheaper than the traditional method on the construction of coarse spaces. The Jacobian matrix in Newton method is not explicitly formed for saving memory, and a preconditioning matrix approximate to the Jacobian matrix is constructed using the collision and streaming terms and ignoring the energy coupling and the angle coupling. By exploring the structure of the matrix in the algebraic transport system, we find that the structures of the submatrices in the preconditioning matrix are similar to each other. Taking advantage of the similarity, we construct subinterpolations based on a submatrix using existing coarsening algorithms, e.g., GAMG in PETSc [13], BoomerAMG in HYPRE [14], etc., for saving compute time. The subinterpolations are then extended to cover the full matrix, and the coarse spaces are finally formed in a Galerkin manner using the extended interpolations. We numerically show that the proposed parallel transport method is scalable to more than 10,000 processor cores for 3D unstructured meshes problems with billions of unknowns. Compared with the traditional multilevel methods in which coarse spaces are built by coarsening the full matrix, the new multilevel methods has a much smaller preconditioner setup time.

The rest of this paper is organized as follows. In Section 2, the multigroup neutron transport equations and its spatial and angular discretizations are presented. And a highly parallel multilevel domain decomposition preconditioner together with a subspace-based coarsening algorithm is discussed in Section 3. In Section 4, some numerical tests are carefully studied to demonstrate the performance of the proposed algorithm. A few remarks and conclusions are drawn in Section 6.

## 2. PROBLEM DESCRIPTION

In this Section, we describe the multigroup neutron transport equations and its stabilizing technique, where the equations is discretized using the first-order Lagrange finite element method in space and using the discrete ordinates scheme in angle.

### 2.1. Multigroup neutron transport equations

As stated earlier, the multigroup neutron transport equations describes the interactions of neutrons with the background materials. The fundamental quantities of interest in the neutron transport calculations are angular fluxes, $\Psi_g$ [ $\text{cm}^{-2}\,\text{s}^{-1}\,\text{st}^{-1}$], $g = 1, 2, ..., G$, on $\mathcal{D} \times \mathcal{S}$, where $G$ is the number of energy groups, $\mathcal{D}$ is a 3D spatial domain as shown in Fig. 1, and $\mathcal{S}$ is a 2D unit sphere of neutron flying directions. The angular fluxes are obtained by solving the following multigroup neutron transport equations

$$\begin{aligned}\vec{\Omega} \cdot \vec{\nabla} \Psi_g + \Sigma_{t,g} \Psi_g &= \sum_{g'=1}^{G} \int_{\mathcal{S}} \Sigma_{s,g' \to g} f_{g' \to g}(\vec{\Omega}' \cdot \vec{\Omega}) \Psi_{g'}(\vec{x}, \vec{\Omega}') \, d\Omega' \\ &\quad + \frac{1}{4\pi} \frac{\chi_g}{k} \sum_{g'=1}^{G} \nu \Sigma_{f,g'} \Phi_{g'}, \text{ in } \mathcal{D} \times \mathcal{S}.\end{aligned} \quad (1)$$

Here $\vec{\Omega} \in \mathcal{S}$ is the independent angular variable for the neutron flying direction, $\vec{x} \in \mathcal{D}$ [cm] represents the independent spatial variable, $\Sigma_{t,g}$ [$\text{cm}^{-1}$] is the macroscopic total cross section, $\Sigma_{s,g' \to g}$ [$\text{cm}^{-1}$] is the macroscopic scattering cross section from group $g'$ to group $g$, $\Sigma_{f,g}$ [$\text{cm}^{-1}$] is the macroscopic fission cross section, $\chi_g$ is the prompt fission spectrum, and $\nu$ is the averaged number of neutrons emitted per fission. $\Phi_g$ [ $\text{cm}^{-2}\,\text{s}^{-1}$] is the scalar flux defined as

$$\Phi_g \equiv \int_{\mathcal{S}} \Psi_g \, d\Omega.$$

$f_{g' \to g}$ [$\text{st}^{-1}$] is the scattering phase function for redistributing neutrons from the incoming directions $\vec{\Omega}'$ to a certain outgoing direction $\vec{\Omega}$, and satisfies

$$\int_{\mathcal{S}} f_{g' \to g}(\vec{\Omega}, \vec{\Omega}') \, d\Omega = 1.$$

$k$ is the eigenvalue (the largest eigenvalue is referred to as the neutron multiplication factor). The largest eigenvalue and its corresponding eigenvector are often referred to as the fundamental mode. Eq. (1) is known as the *k*-eigenvalue problem. The *k*-eigenvalue problem is used to provide the initial condition for fast transient simulations and a factor of $\frac{1}{k}$ is applied to the fission cross section to make sure the initial condition is self-sustained without changing the problem configuration during the transient. In Eq. (1), the first term is the *streaming term*, and the second is the *collision* term. The first term on the right hand side is the *scattering term*, which couples the angular fluxes of all directions and energy groups together. The second term on the right hand side is the fission term coupling the angular fluxes of all directions.

For simplifying discussion, we define some notations that will be used in the rest of the paper.

$$\mathbf{\Psi} \equiv \left[\Psi_1, \Psi_2, \cdots, \Psi_G\right]^T,$$

$$\mathbf{L}\mathbf{\Psi} \equiv \left[\mathbb{L}_1\Psi_1, \mathbb{L}_2\Psi_2, \cdots, \mathbb{L}_G\Psi_G\right]^T, \quad \mathbb{L}_g\Psi_g \equiv \vec{\nabla} \cdot \vec{\Omega}\Psi_g + \Sigma_{t,g}\Psi_g,$$

$$\mathbf{S}\mathbf{\Psi} \equiv \left[\mathbb{S}_1\Psi_1, \mathbb{S}_2\Psi_2, \cdots, \mathbb{S}_G\Psi_G\right]^T, \quad \mathbb{S}_g\Psi_g \equiv \sum_{g'=1}^{G} \int_{\mathcal{S}} \Sigma_{s,g'\to g} f_{g'\to g} \Psi_{g'} \, d\Omega',$$

$$\mathbf{F}\mathbf{\Psi} \equiv \left[\mathbb{F}_1\Psi_1, \mathbb{F}_2\Psi_2, \cdots, \mathbb{F}_G\Psi_G\right]^T, \quad \mathbb{F}_g\Psi_g \equiv \frac{1}{4\pi}\chi_g \sum_{g'=1}^{G} \nu\Sigma_{f,g'}\Phi_{g'}.$$

Here $\mathbf{L}$ is the streaming-collision operator, $\mathbf{S}$ is the scattering operator and $\mathbf{F}$ is the fission operator. With these notations, (1) is rewritten as

$$\mathbf{L}\mathbf{\Psi} = \mathbf{S}\mathbf{\Psi} + \frac{1}{k}\mathbf{F}\mathbf{\Psi}, \tag{2}$$

*2.2. Stabilization techniques*

To present the weak form of Eq. (2), we define a function inner product over $\mathcal{D} \times \mathcal{S}$

$$(\mathbf{a}, \mathbf{b})_{\mathcal{D}\times\mathcal{S}} \equiv \sum_{g=1}^{G} \int_{\mathcal{S}} d\Omega \int_{\mathcal{D}} dx \, a_g(\vec{x}, \vec{\Omega}) b_g(\vec{x}, \vec{\Omega}),$$

where $\mathbf{a} = [a_1, a_2, ..., a_G]^T$ and $\mathbf{b} = [b_1, b_2, ..., b_G]^T$ are generic multigroup functions on $\mathcal{D} \times \mathcal{S}$. The subscript $\mathcal{D} \times \mathcal{S}$ is dropped for simplicity. Similarly, we define the boundary integral as:

$$\langle \mathbf{a}, \mathbf{b} \rangle \equiv \langle \mathbf{a}, \mathbf{b} \rangle^+ + \langle \mathbf{a}, \mathbf{b} \rangle^-, \quad \langle \mathbf{a}, \mathbf{b} \rangle^\pm \equiv \sum_{g=1}^{G} \oint_{\partial\mathcal{D}} dx \int_{\mathcal{S}^\pm_{\vec{n}_b}} d\Omega \left|\vec{\Omega}\cdot\vec{n}_b\right| a_g(\vec{x},\vec{\Omega}) b_g(\vec{x},\vec{\Omega}).$$

Here $\partial\mathcal{D}$ is the boundary of $\mathcal{D}$, $\vec{n}_b$ is the outward unit normal vector on the boundary, and $\mathcal{S}^\pm_{\vec{n}_b} = \{\vec{\Omega} \in \mathcal{S} : \vec{\Omega}\cdot\vec{n}_b \gtrless 0\}$. Following the standard finite element technique, the weak form of Eq. (2) is written as

$$(\mathbf{L}^*\mathbf{\Psi}^*, \mathbf{\Psi}) + \langle \mathbf{\Psi}^*, \bar{\mathbf{\Psi}} \rangle^+ - \langle \mathbf{\Psi}^*, \bar{\mathbf{\Psi}} \rangle^- = (\mathbf{\Psi}^*, \mathbf{S}\mathbf{\Psi}) + \frac{1}{k}(\mathbf{\Psi}^*, \mathbf{F}\mathbf{\Psi}), \tag{3}$$

where $\mathbf{L}^*$ is the adjoint operator of $\mathbf{L}$, that is,

$$\mathbf{L}^*\mathbf{\Psi} \equiv \left[\mathbb{L}^*_1\Psi_1, \mathbb{L}^*_2\Psi_2, \cdots, \mathbb{L}^*_G\Psi_G\right]^T, \quad \mathbb{L}^*_g\Psi_g \equiv -\vec{\nabla}\cdot\vec{\Omega}\Psi_g + \Sigma_{t,g}\Psi_g,$$

and $\mathbf{\Psi}^*$ is a test function. For the boundary terms $\langle \mathbf{\Psi}^*, \bar{\mathbf{\Psi}} \rangle^+ - \langle \mathbf{\Psi}^*, \bar{\mathbf{\Psi}} \rangle^-$, the vacuum boundary condition and the reflecting boundary condition are employed in this study. In the vacuum boundary condition, $\bar{\mathbf{\Psi}}$ is defined as

$$\bar{\mathbf{\Psi}} = \begin{cases} \mathbf{\Psi}, & \text{on } \partial\mathcal{D}, \ \vec{\Omega}\cdot\vec{n}_b \geq 0, \\ \mathbf{0}, & \text{on } \partial\mathcal{D}, \ \vec{\Omega}\cdot\vec{n}_b < 0. \end{cases}$$

For the reflecting boundary condition, $\bar{\mathbf{\Psi}}$ is written as

$$\bar{\mathbf{\Psi}} = \begin{cases} \mathbf{\Psi}, & \text{on } \partial \mathcal{D},\ \vec{\Omega}\cdot\vec{n}_{\text{b}} \geq 0, \\ \mathbf{\Psi}_r, & \text{on } \partial \mathcal{D},\ \vec{\Omega}\cdot\vec{n}_{\text{b}} < 0, \end{cases}$$

where $\mathbf{\Psi}_r$ is the reflecting angular fluxes of $\mathbf{\Psi}$ on $\vec{\Omega}_r = \vec{\Omega} - 2(\vec{\Omega}\cdot\vec{n}_{\text{b}})\vec{n}_{\text{b}}$. The weak form (3) is usually unstable and a stabilizing technique is required. In this work, a SUPG-like (Streamline upwind/Petrov-Galerkin) technique, SAAF (self-adjoint angular flux), is employed. In the SAAF method, the streaming-collision operator $\mathbf{L}$ is split into two parts (the streaming operator $\mathbf{L}_1$ and the collision operator $\mathbf{L}_2$),

$$\mathbf{L}\mathbf{\Psi} \equiv \mathbf{L}_1\mathbf{\Psi} + \mathbf{L}_2\mathbf{\Psi}, \tag{4}$$

where

$$\mathbf{L}_1\mathbf{\Psi} \equiv \left[\mathbb{L}_{1,1}\Psi_1, \mathbb{L}_{1,2}\Psi_2, \cdots, \mathbb{L}_{1,G}\Psi_G\right]^T,\ \mathbb{L}_{1,g}\Psi_g \equiv \vec{\Omega}\cdot\vec{\nabla}\Psi_g,$$

$$\mathbf{L}_2\mathbf{\Psi} \equiv \left[\mathbb{L}_{2,1}\Psi_1, \mathbb{L}_{2,2}\Psi_2, \cdots, \mathbb{L}_{2,G}\Psi_G\right]^T,\ \mathbb{L}_{2,g}\Psi_g \equiv \Sigma_{\text{t},g}\Psi_g.$$

The "inverse" of $\mathbf{L}_2$ is further defined as

$$\mathbf{L}_2^{-1}\mathbf{\Psi} \equiv \left[\mathbb{L}_{2,1}^{-1}\Psi_1, \mathbb{L}_{2,2}^{-1}\Psi_2, \cdots, \mathbb{L}_{2,G}^{-1}\Psi_G\right]^T,\ \mathbb{L}_{2,g}^{-1}\Psi_g = \Psi_g/\Sigma_{\text{t},g}.$$

The corresponding adjoint operators for $\mathbf{L}_1$ and $\mathbf{L}_2$ are simply written as $\mathbf{L}_1^* = -\mathbf{L}_1$ and $\mathbf{L}_2^* = \mathbf{L}_2$, respectively. With these notations, the stabilized weak form obtained using the SAAF method reads as

$$\mathbb{a}\left(\mathbf{\Psi}^*, \mathbf{\Psi}\right) = \frac{1}{k}\mathbb{f}\left(\mathbf{\Psi}^*, \mathbf{\Psi}\right), \tag{5}$$

with

$$\mathbb{a}\left(\mathbf{\Psi}^*, \mathbf{\Psi}\right) \equiv \left(\mathbf{L}_1\mathbf{\Psi}^*, (\tau\mathbf{L}_1 - \mathbf{I} + \tau\mathbf{L}_2)\mathbf{\Psi}\right) + \left(\mathbf{L}_2\mathbf{\Psi}^*, \mathbf{\Psi}\right) + \langle\mathbf{\Psi}^*, \bar{\mathbf{\Psi}}\rangle^+$$
$$- \left((\mathbf{I}+\tau\mathbf{L}_1)\mathbf{\Psi}^*, \mathbf{S}\mathbf{\Psi}\right) - \langle\mathbf{\Psi}^*, \bar{\mathbf{\Psi}}\rangle^-,$$

$$\mathbb{f}\left(\mathbf{\Psi}^*\right) \equiv \left((\mathbf{I}+\tau\mathbf{L}_1)\mathbf{\Psi}^*, \mathbf{F}\mathbf{\Psi}\right).$$

Here $\tau$ is the stabilization parameter defined as

$$\tau \equiv \begin{bmatrix} \tau_1 \\ \tau_2 \\ \vdots \\ \tau_G \end{bmatrix},$$

where

$$\tau_g = \begin{cases} \frac{1}{c\Sigma_{\text{t},g}}, & ch\Sigma_{\text{t},g} \geq \varsigma, \\ \frac{h}{\varsigma}, & ch\Sigma_{\text{t},g} < \varsigma, \end{cases} \tag{6}$$

where $h$ is the characteristic length of a mesh element; $\varsigma$ is usually chosen to be a constant of 0.5. $c$ is a constant, and it is 1.0 by default. Eq. (5) is discretized using the first order finite element method in space and using the discrete ordinates ($S_N$) scheme (can be thought as a collocation method) in angle. More precisely, with an angular quadrature set, $\left\{\vec{\Omega}_d, w_d, d = 1, \cdots, N_d\right\}$, where $N_d$ is the number of angular directions, the multigroup transport equations are calculated along these directions and all angular integrations are approximated with the angular quadrature. That is, an integral of general functions over $\mathcal{S}$ is written as a weighted summation:

$$\int_{\mathcal{S}} \Psi_g\, d\Omega = \sum_{d=1}^{N_d} w_d \Psi_{g,d}.$$

Applying this technique to Eq. (5), we have the following formula for the collision term

$$(\mathbf{L}_2\mathbf{\Psi}^*, \mathbf{\Psi}) = \sum_{g=1}^{G} \int_{\mathcal{S}} d\Omega \left(\Sigma_{t,g}\Psi_g^*, \Psi_g\right)_{\mathcal{D}} = \sum_{g=1}^{G} \sum_{d=1}^{N_d} w_d \left(\Sigma_{t,g}\Psi_{g,d}^*, \Psi_{g,d}\right)_{\mathcal{D}}, \quad (7)$$

where $(\cdot,\cdot)_{\mathcal{D}}$ denotes the integral over $\mathcal{D}$. Similarly, the $S_N$ method can be applied to the fission, the streaming and the scattering terms, and the resulting system has $N_d \times G$ spatial functions on $\mathcal{D}$. We ignore the description of spatial discretization since it is a standard technique and can be found in existing literatures [15]. In next Section, we study a scalable parallel solver with a subspace-based coarsening algorithm for solving the discretized version of Eq. (5).

## 3. PARALLEL SOLVER WITH A SUBSPACE-BASED COARSENING ALGORITHM

In this Section, we describe a parallel Jacobian-free Newton-Krylov method [16] for the eigenvalue calculations of Eq. (5), where Newton method is employed for solving a nonlinear system of equations arising from the reformation of the eigenvalue problems, and during each Newton iteration, a Krylov subspace method such as GMRES [12, 17] is used for calculating the Jacobian system. In order to maintain the scalability of parallel calculations, we introduce some novel coarse spaces that are constructed using a new subspace-based coarsening algorithm for the multilevel domain decomposition preconditioners.

### 3.1. Newton-based eigenvalue solver

The discretized version of Eq. (5), a generalized eigenvalue problem, is written as

$$\mathcal{A}\mathbf{\Psi} = \frac{1}{k}\mathcal{B}\mathbf{\Psi}. \quad (8)$$

Here, without introducing any confusion, the same notation, $\mathbf{\Psi}$, is employed to represent the discretized version of neutron angular fluxes defined on a meh $\mathcal{D}_h$. $\mathcal{A}$ is a matrix arising from the discretization of $\mathtt{a}$ in both space and angle, and $\mathcal{B}$ corresponds to $\mathtt{f}$. One of the simplest algorithms for computing the largest eigenvalue of Eq. 8 is inverse power iteration [18]. The inverse power iteration starts with an initial pair, $(\mathbf{\Psi}_0, k_0 = \|\mathcal{B}\mathbf{\Psi}_0\|)$, and scales the right hand size in place, $\mathcal{B}\mathbf{\Psi}_0 \leftarrow \frac{1}{k_0}\mathcal{B}\mathbf{\Psi}_0$, a new pair, $(\mathbf{\Psi}_{n+1}, k_{n+1})$, is found as follows:

$$\mathcal{A}\mathbf{\Psi}_{n+1} = \mathcal{B}\mathbf{\Psi}_n, \quad (9a)$$

$$k_{n+1} = \|\mathcal{B}\mathbf{\Psi}_{n+1}\|, \quad (9b)$$

$$\mathcal{B}\mathbf{\Psi}_{n+1} \leftarrow \frac{1}{k_{n+1}}\mathcal{B}\mathbf{\Psi}_{n+1}. \quad (9c)$$

Here $n = 1, 2, ..., \max_e$, "$\leftarrow$" represents that the corresponding vector is scaled in place, and $\max_e$ is the maximum number of inverse power iterations. The inverse power iteration works well when the largest eigenvalue and the second largest eigenvalue are not close to each other, but the algorithm converges slow for realistic applications since the largest and the second largest eigenvalues are very close to each other. In this situation, a Newton method is employed to accelerate the inverse power iteration. More precisely, Eqs. (9a), (9b) and (9c) are reformed as a single nonlinear system of equations:

$$\mathcal{F}(\mathbf{\Psi}) = \mathcal{A}\mathbf{\Psi} - \frac{1}{\|\mathcal{B}\mathbf{\Psi}\|}\mathcal{B}\mathbf{\Psi}. \quad (10)$$

An inexact Newton [19] is used to efficiently solve Eq. (10). That is, with a given initial guess $\mathbf{\Psi}_n$, a new approximation $\mathbf{\Psi}_{n+1}$ is formed as follows:

$$\mathbf{\Psi}_{n+1} = \mathbf{\Psi}_n + \alpha_n \Delta \mathbf{\Psi}_n. \tag{11}$$

where $\alpha_n$ is the step size usually calculated using line search methods [20], e.g., back tracking, critical point, etc., and $\Delta \mathbf{\Psi}_n$ is obtained by solving the following Jacobian system

$$\mathcal{J}(\mathbf{\Psi}_n) \Delta \mathbf{\Psi}_n = -\mathcal{F}(\mathbf{\Psi}_n). \tag{12}$$

Here $\mathcal{F}(\mathbf{\Psi}_n)$ is the nonlinear function evaluated at $\mathbf{\Psi}_n$, and $\mathcal{J}(\mathbf{\Psi}_n)$ is the Jacobian matrix at $\mathbf{\Psi}_n$. As stated earlier, there are $N_d \times G$ spatial functions on $\mathcal{D}$ after the angular discretization, where $G$ is the number of energy groups and $N_d$ is the number of angular quadrature points. If $\mathcal{J}$ is formed explicitly, there will be a $(N_d G) \times (N_d G)$ dense block matrix for each mesh vertex and therefore a significant amount of memory is required to store $\mathcal{J}$. More important, it is nontrivial to compute the explicit form of the full matrix of $\mathcal{J}$ since the derivatives of $1/\|\mathcal{B}\mathbf{\Psi}\|$ are challenging to compute. To overcome this difficulty, we employ the Jacobian-free version of the inexact Newton [21]. When solving the Jacobian system (12), GMRES is used and during each GMRES iteration, a matrix-vector multiplication, $\mathcal{J}(\mathbf{\Psi}_n) \Delta \mathbf{\Psi}_n$, is required. In the Jacobian-free method, this matrix-vector multiplication is carried out as

$$\mathcal{J}(\mathbf{\Psi}_n) \Delta \mathbf{\Psi}_n = \frac{\mathcal{F}(\mathbf{\Psi}_n + \delta \Delta \mathbf{\Psi}_n) - \mathcal{F}(\mathbf{\Psi}_n)}{\delta},$$

where $\delta$ is a small permutation, and each matrix-vector multiplication involves one nonlinear function evaluation, that is, $\mathcal{F}(\mathbf{\Psi}_n + \delta \Delta \mathbf{\Psi}_n)$. It is well-known that GMRES convergence is slow when the condition number of $\mathcal{J}$ is high, which is true for most of realistic simulations. Next, we will discuss an efficient and scalable parallel preconditioner based on domain decomposition methods to fix this issue.

*3.2. Multilevel preconditioner with a subspace-based coarsening algorithm*

Below, we describe a domain decomposition framework for constructing a scalable parallel preconditioner, and then introduce a subspace-based coarsening algorithm for building coarse spaces. To apply a preconditioner, we rewrite Eq. (12) as its preconditioned form

$$\mathcal{J}\mathcal{P}^{-1}\mathcal{P}\Delta\mathbf{\Psi} = -\mathcal{F}, \tag{13}$$

where $\mathcal{P}$ is the preconditioning matrix that is an approximation of $\mathcal{J}$, and $\mathcal{P}^{-1}$ is a preconditioning process. Eq. (13) is carried out in two steps, that is,

$$\mathcal{J}\mathcal{P}^{-1}\Delta\tilde{\mathbf{\Psi}} = -\mathcal{F}, \tag{14a}$$

$$\mathcal{P}\Delta\mathbf{\Psi} = \Delta\tilde{\mathbf{\Psi}}, \tag{14b}$$

where $\Delta\tilde{\mathbf{\Psi}}$ is an auxiliary vector. In Eq. (13), the preconditioner is applied from the right, but we want to mention that the preconditioner can be applied from the left as well. General speaking, a preconditioning process is designed to find the solution of the following residual equations,

$$\mathcal{P}\mathbf{e} = \mathbf{r}, \tag{15}$$

where $\mathbf{e}$ is a correction vector, and $\mathbf{r}$ is the residual vector from the outer solver GMRES. In this paper, Eq. (15) is implemented based on domain decomposition methods. We denote a triangulation of $\mathcal{D}$ as $\mathcal{D}_h$. The basic idea of the domain decomposition methods is to partition $\mathcal{D}_h$ into $np$ submeshes $\mathcal{D}_{h,i}, i = 1, 2, ..., np$, ($np$ is the number of processor cores), $np$ submesh problems are solved independently in parallel, and a global correction is formed by combining all the submesh solutions together. To partition $\mathcal{D}_h$, a hierarchal partitioning algorithm [22, 23] is employed since the existing partitioners such as ParMETIS [24] are from ideal when the number of processor cores is large. The hierarchal partitioning takes into consideration that there are multiple processor cores on each compute node in modern supercomputers and all the in-node cores share the memory. The hierarchical partitioning is implemented by applying an existing partitioner such as ParMETIS or PT-Scotch [25] twice or multiple times. More precisely, $\mathcal{D}_h$ is partitioned into $np_1$ submeshes ($np_1$ is the number of compute nodes), and each submesh is further divided into $np_2$ smaller submeshes ($np_2$ is the number of processor cores per compute node). The basic idea of the hierarchal partitioning is quite simple, but it is very effective since it reduces the communication between the compute nodes and improves the partition quality. A mesh sample partitioned into 40 submeshes is shown in Fig. 1, where we assume there are 4 compute nodes, and each node has 10 processor cores. The mesh is partitioned into 4 submeshes, each submesh is further split into 10 small submeshes, and 40 small submeshes in total are obtained. For parallel simulations, the global matrix and the global vector are constructed in such a way that their partitions are consistent with the mesh partition. The unknowns are defined based on mesh vertex with the first order finite element method, while the mesh partition is carried out based on mesh element. After the mesh partition, some mesh vertices are shared by multiple processor cores, and the unknowns associated with these shared vertices are allocated to the processor core with the smallest MPI rank in libMesh [26]. This approach works well when the number of unknowns is small, but it leads to an imbalanced calculation when the number of unknowns is large. To fix this issue, we employ a partition-based assigning scheme, where a partitioner is applied to partition the shared surface mesh into two parts, and each part is then assigned to one processor core. The partition-based scheme is able to generate balanced calculations meanwhile preserving the data locality that is critical for the preconditioner performance. Interested readers are referenced to [23, 27] for a detailed description of the algorithm. To define a domain decomposition-based preconditioner, we introduce some notations. Nonoverlapping submatrices and subvectors on $\mathcal{D}_{h,i}$ are denoted as $\mathcal{P}_i$, $\mathbf{e}_i$ and $\mathbf{r}_i$, respectively. Local subvectors $\mathbf{r}_i$ are extended to overlap with its neighbors by $\delta$ layers according to the nonzero pattern of $\mathcal{P}$, and the overlapping subvectors are denoted as $\mathbf{r}_i^\delta$. A restriction operator extracting overlapping local components from the global vector is defined as

$$\mathbf{r}_i^\delta = \mathcal{R}_i^\delta \mathbf{r} = \begin{pmatrix} I & 0 \end{pmatrix} \begin{pmatrix} \mathbf{r}_i^\delta \\ \mathbf{r}/\mathbf{r}_i^\delta \end{pmatrix},$$

where $I$ is an identity matrix having the same size as $\mathbf{r}_i^\delta$, and $\mathbf{r}/\mathbf{r}_i^\delta$ represents all the components in $\mathbf{r}$ but not in $\mathbf{r}_i^\delta$. The nonoverlapping version of $\mathcal{R}_i^\delta$ is denoted as $\mathcal{R}_i^0$. With those notations, the one-level domain decomposition based preconditioner (restricted Schwarz preconditioner [28]) is written as

$$\mathcal{P}_{one}^{-1} = \sum_{i=1}^{np} (\mathcal{R}_i^0)^T (\mathcal{P}_i^\delta)^{-1} \mathcal{R}_i^\delta, \quad \mathcal{P}_i^\delta = \mathcal{R}_i^\delta \mathcal{P} (\mathcal{R}_i^\delta)^T, \tag{16}$$

where $(\mathcal{P}_i^\delta)^{-1}$ is a subdomain solver such as ILU, or SOR. We want to mention that for different variants of domain decomposition methods, interested readers are referred to [29, 30, 31]. The one-level preconditioner works well when the number of processor cores is small, and coarse spaces are usually required for a large number of processor cores. The coarse spaces can be constructed geometrically [22, 32, 33, 27, 34] or algebraically [35, 36]. The construction of coarse spaces is often expensive and not ideally scalable for realistic simulations. To maintain the parallel scalability of the neutron transport calculations, we introduce a novel subspace-based coarsening algorithm to construct coarse spaces.

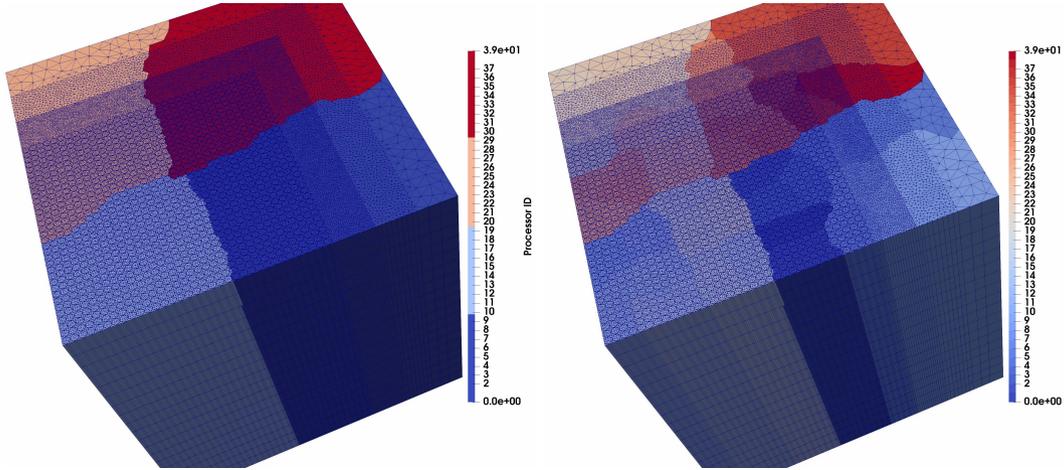

Figure 1. Demonstration of partitioning a mesh into 40 submeshes. Each compute node has 10 processor cores, and 4 compute nodes are available. Each submesh is divided into 10 small submeshes.

To describe the subspace-based coarsening algorithm, we explore the structure of the preconditioning matrix $\mathcal{P}$. For saving memory, $\mathcal{P}$ accounts for only the first and the second terms of Eq (5). This treatment results in a block diagonal matrix, and each block corresponds to the spatial discretization. Instead of coarsening $\mathcal{P}$, the coarsening algorithm is applied to one block and a sequence of subinterpolation matrices are generated. The subinterpolations are extended to cover the full space. Below, we describe this algorithm in detail. If the unknowns were ordered group-by-group and direction-by-direction, $\mathcal{P}$ could be rewritten as

$$\mathcal{P} = \begin{bmatrix} \begin{bmatrix} \mathcal{P}_{1,1}^{(1,1)} & & & \\ & \mathcal{P}_{1,1}^{(2,2)} & & \\ & & \ddots & \\ & & & \mathcal{P}_{1,1}^{(N_d,N_d)} \end{bmatrix} & & \\ & \ddots & \\ & & \begin{bmatrix} \mathcal{P}_{G,G}^{(1,1)} & & & \\ & \mathcal{P}_{G,G}^{(2,2)} & & \\ & & \ddots & \\ & & & \mathcal{P}_{G,G}^{(N_d,N_d)} \end{bmatrix} \end{bmatrix}, \tag{17}$$

where a big block represents the coupling matrix within a given energy group, and $\mathcal{P}_{g,g'}^{(d,d')}$ represents the coupling between angular directions $d$ and $d'$ in the $g$th group. If the scattering and the fission terms were taken into account, $\mathcal{P}$ would be a fully coupled matrix. $\mathcal{P}_{g,g'}^{(d,d')}$ is nonzero only when $g = g'$ and $d = d'$, otherwise, $\mathcal{P}_{g,g'}^{(d,d')} = 0$. To simplify description, the submatrices $\mathcal{P}_{g,g}^{(d,d)}$ are rewritten as $\mathcal{P}_{(j)}$, $j = g \times N_d + d$, $g = 1, 2, .., G$, and $d = 1, 2, .., N_d$. As state earlier, the structures of $\mathcal{P}_{(j)}$ are similar to each other since they come from the same operators, but their numeric values are different from each other because material properties (cross sections) for each energy group are different. The similarity of the submatrix structures implies that the subinterpolations for a given submatrix can be used for other submatrices as well. The difference of the numeric values of the submatrices indicates that coarse spaces for a given submatrix should not be used in other submatrices. A submatrix can be coarsened using existing algorithms such as BoomerAMG in HYPRE [14] and GAMG in PETSc [13]. Let us assume that a $L$-level method is generated using an existing algorithm to coarsen the first submatrix $\mathcal{P}_{(1)}$, and the corresponding $(L-1)$ subinterpolations are denoted as $\tilde{\mathcal{I}}_{l+1}^l, l = 1, 2, ..., L-1$. The $L$th level is the coarsest level, and the 1st level is the finest level. We denote the full vectors on the $l$th level as $\mathbf{r}^l, \mathbf{e}^l, l = 1, 2, ..., L$ ($\mathbf{r}^1 = \mathbf{r}$ and $\mathbf{e}^1 = \mathbf{e}$ ). We define a restriction operator $\mathcal{R}_{(j)}^l$ on the $l$th level, which extracts the $j$th subvector $\mathbf{r}_{(j)}^l$ for the $d$th angular direction of the $g$th energy group from the full vector $\mathbf{r}^l$, that is,

$$\mathbf{r}_{(j)}^l = \mathcal{R}_{(j)}^l \mathbf{r}^l = \begin{pmatrix} I_{(j)}^l & 0 \end{pmatrix} \begin{pmatrix} \mathbf{r}_{(j)}^l \\ \mathbf{r}^l / \mathbf{r}_{(j)}^l \end{pmatrix}, l = 1, 2, .., L$$

where $I_{(j)}^l$ is an identiy matrix that has the same size as $\mathbf{r}_{(j)}^l$, and $\mathbf{r}^l / \mathbf{r}_{(j)}^l$ represents all the components in $\mathbf{r}^l$ but not in $\mathbf{r}_{(j)}^l$. With these notations, the full interpolations $\mathcal{I}_{l+1}^l$ are constructed from the subinterpolations $\tilde{\mathcal{I}}_{l+1}^l$ as follows:

$$\mathcal{I}_{l+1}^l = \sum_{j=1}^{G \times N_d} (\mathcal{R}_{(j)}^l)^T \tilde{\mathcal{I}}_{l+1}^l \mathcal{R}_{(j)}^{l+1}, l = 1, 2, ..., L-1. \tag{18}$$

Coarse operators $\mathcal{P}^l$ are then formed in a Galerkin method with these full interpolations $\mathcal{I}_{l+1}^l$, that is,

$$\mathcal{P}^{l+1} = (\mathcal{I}_{l+1}^l)^T \mathcal{P}^l \mathcal{I}_{l+1}^l, l = 1, 2, ..., L-1. \tag{19}$$

Here the spars matrix triple products are implemented using a memory-efficient method studied in [37]. Finally, a multilevel additive Schwarz preconditioner with the subspace-based coarsening algorithm (abbreviated as "MASM$_{\text{sub}}$") is summarized in Alg. 1, where the choice of submatrix to coarsen is arbitrary. For comparison, we also present the traditional multilevel additive Schwarz preconditioner (referred to as "MASM") in Alg. 2, where the coarsening algorithm is carried out using the full matrix. Solving phases in both algorithms are the same, and the setup phases are different. The setup cost of MASM$_{\text{sub}}$ is much smaller than that of MASM since a small submatrix is employed in the coarsening. In next Section, we will verify this statement by comparing the performance of MASM$_{\text{sub}}$ with that of MASM. We want to mention that the multilevel preconditioner with geometric coarse spaces has been successfully applied to different problems in our previous works, e.g., elasticity problems [22], and fluid-structure interactions [32, 27, 34, 33], and incompressible flow problems [38]. However, geometric coarse spaces are not available in this work, and thus the algebraic coarse space have to be constructed.

## 4. NUMERICAL RESULTS

In this Section, we report the performance of the proposed algorithm in terms of the compute time and the parallel efficiency for a standard benchmark problem, C5G7 3D, as shown in Fig. 1 and 2. The three-dimensional domain, shown in Fig. 2, consist of four assemblies, and each assembly is made up of a $17 \times 17$ lattice of pin cells. The overall dimensions of the domain is $64.26 \times 64.26 \times 64.26$ cm (width $\times$ length $\times$ height), while each assembly is $21.42 \times 21.42 \times 42.84$ cm. On the top of the assemblies, there are some control rods. All of these pins have a 0.54 radius. The different colors of pins correspond to different materials, e.g., UO2, guide tube, fission chamber, MOX 4.3%, MOX 7.0%, MOX 8.7%, and control rod. Interested readers are referred to [39]

---

**Algorithm 1** Multilevel additive Schwarz preconditioner with the subspace-based coarsening algorithm (MASM$_{\text{sub}}$)

---

1: **procedure** PCSETUP($\mathcal{P}$)
2:     Extract the first submatrix $\mathcal{P}_{(1)}$ from $\mathcal{P}$ in parallel
3:     Coarsen $\mathcal{P}_{(1)}$ to generate $(L-1)$ subinterpolations $\tilde{\mathcal{I}}_l^{l+1}$
4:     **for** $l = 1, 2, ..., L-1$ **do**
5:         Construct a full interpolation: $\mathcal{I}_{l+1}^l = \sum_{j=1}^{G \times N_d} (\mathcal{R}_{(j)}^l)^T \tilde{\mathcal{I}}_{l+1}^l \mathcal{R}_{(j)}^{l+1}$
6:     **end for**
7:     **for** $l = 1, 2, ..., L-1$ **do**
8:         Build a full coarse matrix using Galerkin method: $\mathcal{P}^{l+1} = (\mathcal{I}_{l+1}^l)^T \mathcal{P}^l \mathcal{I}_{l+1}^l$
9:     **end for**
10:    Return $\mathcal{P}^l$ and $\mathcal{I}_{l+1}^l$
11: **end procedure**
12: **procedure** PCAPPLY($\mathcal{P}^l, \mathbf{e}^l, \mathbf{r}^l$)
13:     **if** $l = L$ **then**
14:         Solve $\mathcal{P}^L \mathbf{e}^L = \mathbf{r}^L$ with a redundant direct solver on each compute node
15:     **else**
16:         Pre-solve $\mathcal{P}^l \mathbf{e}^l = \mathbf{r}^l$ using an iterative solver preconditioned by $(\mathcal{P}_{one}^l)^{-1}$
17:         Set $\bar{\mathbf{r}}^l = \mathbf{r}^l - \mathcal{P}^l \mathbf{e}^l$
18:         Apply the restriction: $\bar{\mathbf{r}}^{l+1} = (\mathcal{I}_{l+1}^l)^T \bar{\mathbf{r}}^l$
19:         $\mathbf{z}^{(l+1)} = \text{MASM}(\mathcal{P}^{l+1}, \mathbf{z}^{l+1}, \bar{\mathbf{r}}^{l+1})$
20:         Apply the interpolation: $\mathbf{z}^l = \mathcal{I}_{l+1}^l \mathbf{z}^{l+1}$
21:         Correct the solution: $\mathbf{e}^l = \mathbf{e}^l + \mathbf{z}^l$
22:         Post-solve $\mathcal{P}^l \mathbf{e}^l = \mathbf{r}^l$ using an iterative solver preconditioned by $(\mathcal{P}_{one}^l)^{-1}$
23:     **end if**
24:     Return $\mathbf{e}^l$
25: **end procedure**

---

**Algorithm 2** Traditional multilevel additive Schwarz preconditioner (MASM)

---

1: **procedure** PCSETUP($\mathcal{P}$)
2:     Coarsen $\mathcal{P}$ to generate $(L-1)$ interpolations $\mathcal{I}_l^{l+1}$ and $(L-1)$ coarse matrices $\mathcal{P}^l$
3:     Return $\mathcal{P}^l$ and $\mathcal{I}_l^{l+1}$
4: **end procedure**
5: **procedure** PCAPPLY($\mathcal{P}^l, \mathbf{e}^l, \mathbf{r}^l$)
6:     **if** $l = L$ **then**
7:         Solve $\mathcal{P}^L \mathbf{e}^L = \mathbf{r}^L$ with a redundant direct solver on each compute node
8:     **else**
9:         Pre-solve $\mathcal{P}^l \mathbf{e}^l = \mathbf{r}^l$ using an iterative solver preconditioned by $(\mathcal{P}_{one}^l)^{-1}$
10:        Set $\bar{\mathbf{r}}^l = \mathbf{r}^l - \mathcal{P}^l \mathbf{e}^l$
11:        Apply the restriction: $\bar{\mathbf{r}}^{l+1} = (\mathcal{I}_{l+1}^l)^T \bar{\mathbf{r}}^l$
12:        $\mathbf{z}^{(l+1)} = \text{MASM}(\mathcal{P}^{l+1}, \mathbf{z}^{l+1}, \bar{\mathbf{r}}^{l+1})$
13:        Apply the interpolation: $\mathbf{z}^l = \mathcal{I}_{l+1}^l \mathbf{z}^{l+1}$
14:        Correct the solution: $\mathbf{e}^l = \mathbf{e}^l + \mathbf{z}^l$
15:        Post-solve $\mathcal{P}^l \mathbf{e}^l = \mathbf{r}^l$ using an iterative solver preconditioned by $(\mathcal{P}_{one}^l)^{-1}$
16:     **end if**
17:     Return $\mathbf{e}^l$
18: **end procedure**

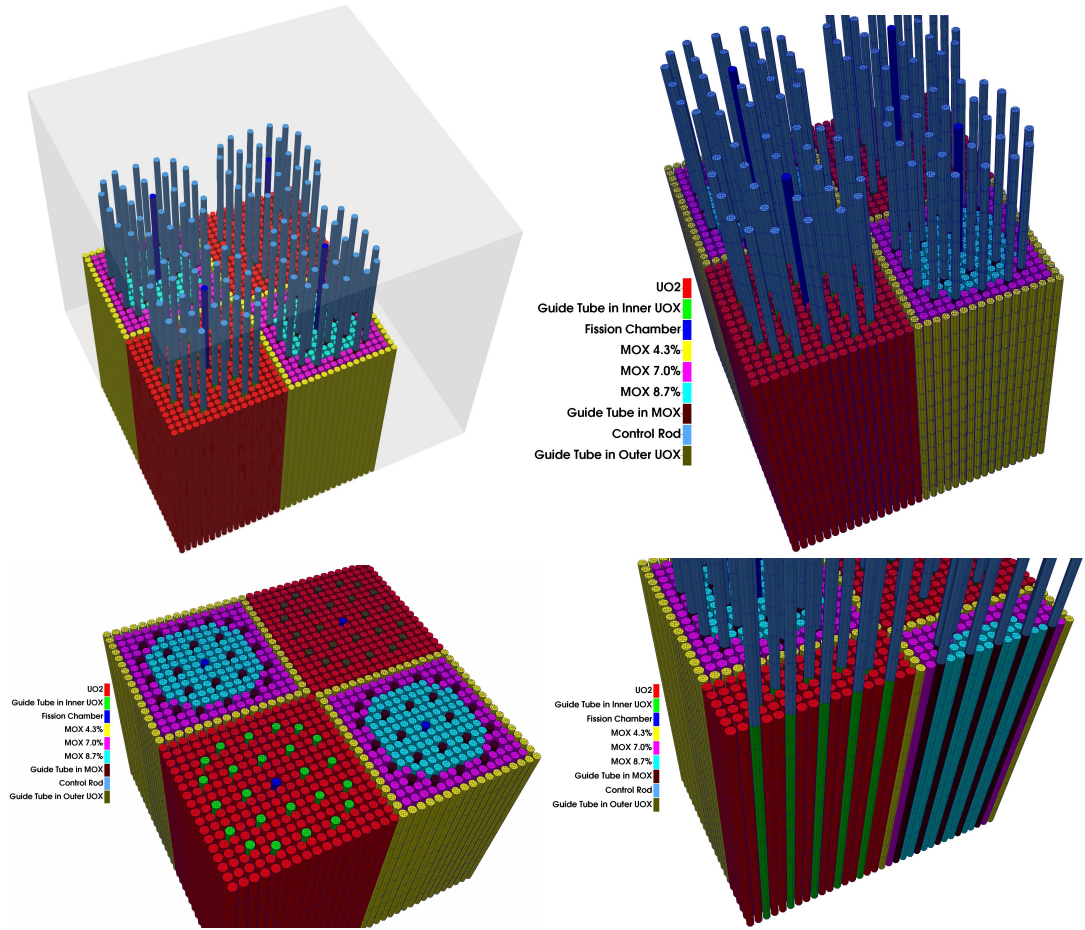

Figure 2. Computational domain of C5G7 3D benchmark. Top left: four assemblies in the front bottom corner of $64.26 \times 64.26 \times 64.26$ cm cube; top right: four assemblies and control rods; bottom left: XY plane cross section of the assemblies; bottom right: XZ plane cross section of the assemblies. Different colors of pins represent different materials.

for more details. Material properties (cross sections) can be found in [39]. The vacuum boundary conditions are applied to the back, the right and the top surfaces, while the reflected boundary conditions are applied to the front, the left and the bottom surfaces. The numerical experiments in this section are carried out on a supercomputer at INL (Idaho National Laboratory), where each compute node has two 20-core processors with 2.4 GHz and the compute nodes are connected by an OmniPath network. The eigenvalue problem (1) is solved on the domain shown in Fig. 2, and the zero order flux moments for the first and the seventh energy groups are shown in Fig. 3, and the computed eigenvalue is 1.1416730884. The discretization of the multigroup neutron transport equations is implmented using RattleSnake [15], MOOSE [40] and libMesh [26]. The proposed preconditioner is implemented in PETSc [13] as part of this work, and the submatrix is coarsened using BoomerAMG in HYPRE[41].

For convenience, we define some notations that will be used in the rest of the paper. "PC" represents the preconditioning algorithms (MASM or $\text{MASM}_{\text{sub}}$). "$\text{Iter}_{\text{Newton}}$" is the number of Newton iterations, and "$\text{Iter}_{\text{GMRES}}$" is the averaged number of GMRES iterations per Newton step. "$\text{Time}_{\text{PCSetup}}$" is the time spent on the preconditioner setup, "$\text{Time}_{\text{PCApply}}$" is the compute time on the application of preconditioner, "$\text{Time}_{\text{KSP}}$" is the compute time spent on the linear solver, and "$\text{Time}_{\text{Total}}$" is the total compute time of the overall simulation. "EFF" is the parallel efficiency. We want to mention that both "Time$_{\text{PCSetup}}$" and "Time$_{\text{PCApply}}$" are part of "Time$_{\text{KSP}}$", and "Time$_{\text{KSP}}$" is part of "Time$_{\text{Total}}$". Therefore, we will see that the times on the linear solver and the overall simulation are reduced when we have an improvement on the preconditioner setup. "Time$_{\text{Func}}$" is the time spent on the function evaluations, "Time$_{\text{Jac}}$" is the time spent on the Jacobian evaluations, "Time$_{\text{LS}}$" is the time spent on the line search method, and "Time$_{\text{MF}}$" is the compute time on matrix-free operations. SOR is used as the submesh solver on each level, the submesh overlapping size is set to zero, a relative tolerance of $10^{-6}$ is chosen for Newton iteration and the Jacobian system is solved inexactly with a relative tolerance of $10^{-1}$.

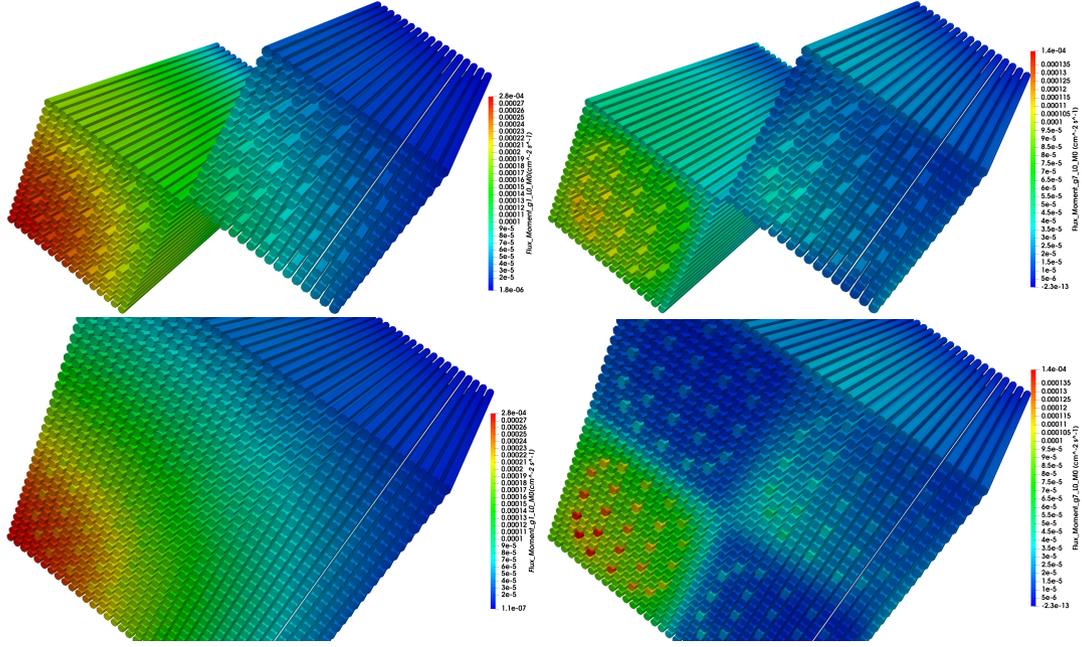

Figure 3. Zero order flux moments. Top left: the zero order flux moment of the first energy group for UO2; To right: the zero order flux moment of the seventh energy group for UO2; bottom left: the zero order flux moment of the first energy group for all the four assemblies; bottom right: the zero order flux moment of the seventh energy group for all the four assemblies.

### 4.1. Performance comparison between MASM and MASM$_{sub}$

In this test, we study the performance of MASM$_{sub}$, compared with the traditional MASM, using 160, 320, 640 and 1,280 processor cores. A mesh with 832,371 nodes and 1,567,944 elements is used. A Gauss Chebyshev quadrature is chosen with 32 angular directions, and there are 224 unknowns per mesh vertex. The resulting system of eigenvalue equations has 186,451,104 unknowns. Two inverse power iterations are used to generate an initial guess for the Newton-Krylov method. Numerical results on Newton iteration, GMRES iteration and compute times of different components of the algorithms are summarized in Table I. We observed, from Table I, that the numbers of iterations for both Newton and GMRES stay as constants when we increase the number of processor cores from 160 to 320, 640 and 1,280. This indicates that the overall algorithm is mathematically scalable. The preconditioner setup time of MASM$_{sub}$ is smaller than that used in MASM for all processor counts. For example, the compute time on the preconditioner setup for MASM$_{sub}$ using 160 processor cores is 33 s only, while that of MASM is 287 s. That is, MASM$_{sub}$ is nine times faster than MASM when we use 160 processor cores. The similar behaviors are observed for all core counts. Preconditioner application times for MASM$_{sub}$ and MASM are similar to each other, and sometimes, MASM$_{sub}$ is slightly faster than the traditional MASM. The compute times on the linear solver and the entire simulations are accordingly decreased due to the time reduction on the preconditioner setup. The new MASM$_{sub}$ improves the parallel scalability for all core counts, compared with the traditional MASM. The total compute time for MASM$_{sub}$ is almost halved, e.g., it is reduced to 1144 s from 2162 s, when the number of processor cores is increased from 160 to 320. It continues being decreased to 609 s and 347 s when we use 640 and 1,280 processor cores. The parallel efficiency of the new proposed algorithm is better than that obtained using the traditional MASM by 20 percentage points when we use up to 1,280 processor cores, e.g., the new MASM$_{sub}$ has a parallel efficiency of 79% at 1,280 cores, while MASM has a 60% parallel efficiency only. The parallel efficiencies and speedups are plotted in Fig. 4 as well. To have an intuitive analysis, in Fig. 5, we summarize a performance comparison between MASM$_{sub}$ and MASM on the preconditioner time and the total compute time for all core counts. It is easily found that MASM$_{sub}$ outperforms the traditional MASM that coarsens the full matrix.

Table I. Performance comparison between MASM and MASM$_{sub}$ using 160, 320, 640 and 1,280 processor cores for a problem with 186,451,104 unknowns. The system of eigenvalue equations is solved using an inexact Newton-Krylov method.

| $np$ | PC | Iter$_{Newton}$ | Iter$_{GMRES}$ | Time$_{PCSetup}$ | Time$_{PCApply}$ | Time$_{KSP}$ | Time$_{Total}$ | EFF |
|---|---|---|---|---|---|---|---|---|
| 160 | MASM | 5 | 20 | 287 | 225 | 2448 | 2756 | – |
| 160 | MASM$_{sub}$ | 5 | 20 | 33 | 224 | 2008 | 2162 | 100% |
| 320 | MASM | 5 | 20 | 124 | 113 | 1261 | 1421 | 76% |
| 320 | MASM$_{sub}$ | 5 | 20 | 19 | 112 | 1064 | 1144 | 95% |
| 640 | MASM | 5 | 20 | 76 | 65 | 691 | 776 | 70% |
| 640 | MASM$_{sub}$ | 5 | 20 | 12 | 56 | 566 | 609 | 89% |
| 1,280 | MASM | 5 | 20 | 60 | 41 | 406 | 452 | 60% |
| 1,280 | MASM$_{sub}$ | 5 | 20 | 10 | 30 | 323 | 347 | 78% |

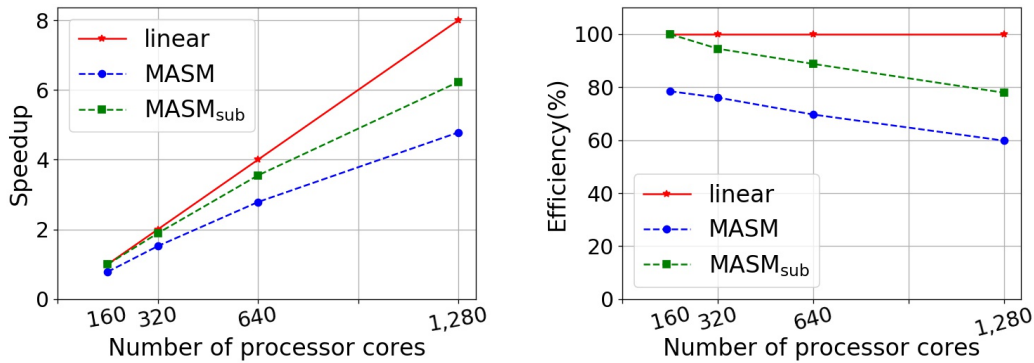

Figure 4. Speedups and parallel efficiencies for MASM and MASM$_{sub}$ on a problem with 186,451,104 unknowns using 160, 320, 640, and 1,280 processor cores, respectively. Left: speedups, right: parallel efficiencies.

### 4.2. Strong scaling study with different problem sizes

In this Section, we will study the performance of the proposed algorithm with respect to different problem sizes on 1,280, 2,560, 5,120 and 10,240 processor cores. The Level-Symmetric angular quadrature is chosen with different numbers of angular directions. A mesh with 6,464,825 vertices and 12,543,552 elements is employed.

We start with a test having 8 angular directions, and the resulting system has 56 variables per mesh vertex. The total number of unknowns in the system of eigenvalue equations is

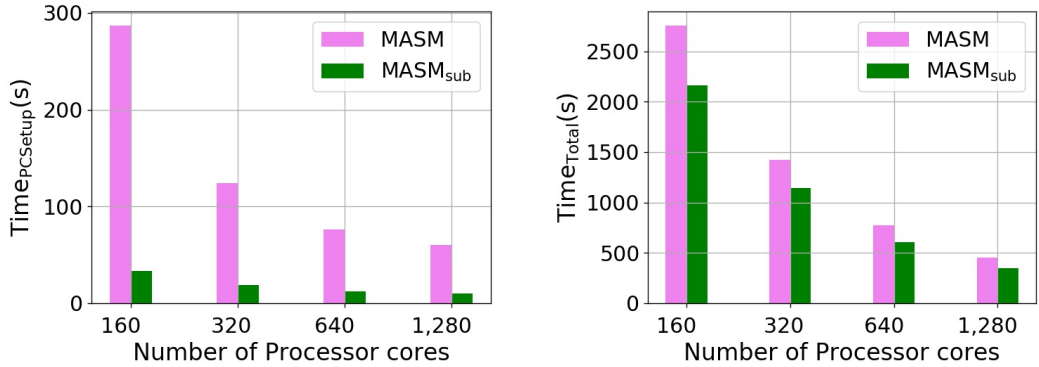

Figure 5. Comparison of preconditioner setup time and total compute time between $\text{MASM}_{\text{sub}}$ and MASM using 160, 320, 640 and 1,280 processor cores for a problem with 186,451,104 unknowns. Left: preconditioner setup time, right: total compute time.

Table II. Performance comparison between $\text{MASM}_{\text{sub}}$ and MASM on a problem with 362,030,200 unknowns using 1,280, 2,560, 5,120 and 10,240 processor cores. A system of eigenvalue equations is computed using an inexact Newton-Krylov method, where $\text{MASM}_{\text{sub}}$ and MASM are employed as a preconditioner, respectively.

| $np$ | PC | $\text{Iter}_{\text{Newton}}$ | $\text{Iter}_{\text{GMRES}}$ | $\text{Time}_{\text{PCSetup}}$ | $\text{Time}_{\text{PCApply}}$ | $\text{Time}_{\text{KSP}}$ | $\text{Time}_{\text{Total}}$ | EFF |
|---|---|---|---|---|---|---|---|---|
| 1,280 | MASM | 8 | 16 | 56 | 69 | 539 | 585 | – |
| 1,280 | $\text{MASM}_{\text{sub}}$ | 7 | 16 | 13 | 55 | 442 | 485 | 100% |
| 2,560 | MASM | 8 | 16 | 40 | 43 | 311 | 336 | 72% |
| 2,560 | $\text{MASM}_{\text{sub}}$ | 7 | 16 | 9 | 30 | 246 | 269 | 90% |
| 5,120 | MASM | 8 | 16 | 40 | 32 | 204 | 218 | 56% |
| 5,120 | $\text{MASM}_{\text{sub}}$ | 7 | 16 | 8 | 18 | 143 | 155 | 78% |
| 10,240 | MASM | 7 | 14 | 37 | 25 | 131 | 140 | 43% |
| 10,240 | $\text{MASM}_{\text{sub}}$ | 7 | 16 | 10 | 32 | 121 | 129 | 47% |

362,030,200 unknowns. The performance comparison between $\text{MASM}_{\text{sub}}$ and MASM for this configuration is shown in Table II, where we observed that the number of Newton iterations is kept close to a constant and the number of GMRES iterations stays close to a constant as well when the number of processor cores is increased from 1,280 to 2,560, 5,120 and 10,240. This trend indicates, again, that both $\text{MASM}_{\text{sub}}$ and MASM make the overall algorithm mathematically scalable. The number of Newton iterations used in $\text{MASM}_{\text{sub}}$ is slightly smaller than that in MASM, and the numbers of GMRES iterations are almost the same for both. We are able to reduce the preconditioner setup time when using $\text{MASM}_{\text{sub}}$ instead of MASM. The time spent on the preconditioner setup for MASM is 56 s at 1,280 cores, and it is reduced to 13 s when we replace MASM with $\text{MASM}_{\text{sub}}$. It is found that this similar behavior happens for other core counts as well, that is, $\text{MASM}_{\text{sub}}$ is four times faster than MASM on the preconditioner setup for all core counts. The preconditioner applicaiton time of $\text{MASM}_{\text{sub}}$ is smaller than MASM for all core counts except $10,240$, where $\text{MASM}_{\text{sub}}$ costs more than MASM by 7 s since two more GMRES iterations are required in $\text{MASM}_{\text{sub}}$. The linear solver time is accordingly decreased due to a time reduction on the preconditioner setup. A good scalability is maintained for $\text{MASM}_{\text{sub}}$ with up to 5,120 processor cores, while that using MASM is not ideally. The problem size is too small for $10,240$ processor cores, and the parallel scalability will be gained if we increase the number of angular directions to have more unknowns (we will test this shortly). The comparison between $\text{MASM}_{\text{sub}}$ and MASM in the preconditioner setup time and the total compute time is also drawn in Fig. 6, where we observed that $MASM_{sub}$ is always better than MASM for all core counts.

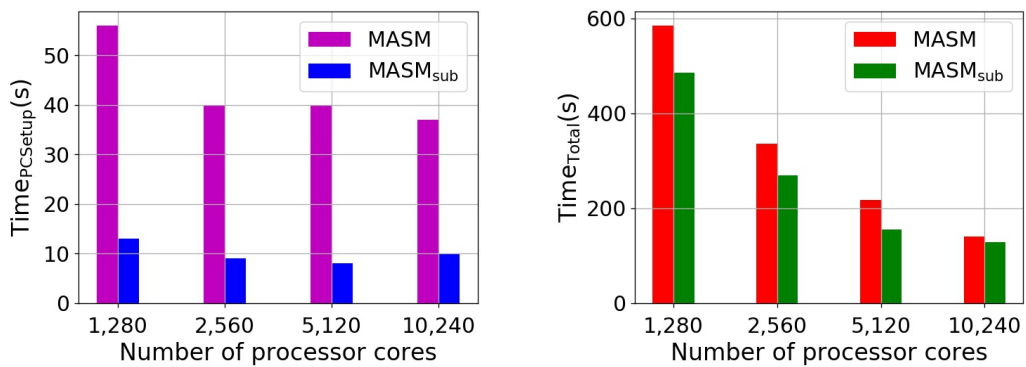

Figure 6. Times on the preconditioner setup and the overall simulation for $MASM_{sub}$ and MASM on a problem with 362,030,200 unknowns using 1,280, 2,560, 5,120 and 10,240 processor cores. Left: preconditioner setup, right: total compute time.

Table III. Performance comparison using 24 angular directions between $MASM_{sub}$ and MASM on 1,280, 2,560, 5,120 and 10,240 processor cores. A system of eigenvalue equations with 1,086,090,600 unknowns is calculated using an inexact Newton-Krylov method.

| $np$ | PC | $Iter_{Newton}$ | $Iter_{GMRES}$ | $Time_{PCSetup}$ | $Time_{PCApply}$ | $Time_{KSP}$ | $Time_{Total}$ | EFF |
|---|---|---|---|---|---|---|---|---|
| 1,280 | MASM | 7 | 14 | 327 | 184 | 1603 | 1740 | – |
| 1,280 | $MASM_{sub}$ | 7 | 16 | 32 | 161 | 1350 | 1485 | 100% |
| 2,560 | MASM | 7 | 14 | 193 | 115 | 881 | 953 | 78% |
| 2,560 | $MASM_{sub}$ | 7 | 16 | 22 | 84 | 720 | 792 | 94% |
| 5,120 | MASM | 7 | 14 | 185 | 107 | 607 | 647 | 57% |
| 5,120 | $MASM_{sub}$ | 7 | 16 | 16 | 46 | 403 | 443 | 84% |
| 10,240 | MASM | 7 | 14 | 168 | 85 | 450 | 474 | 39% |
| 10,240 | $MASM_{sub}$ | 7 | 16 | 19 | 45 | 268 | 292 | 64% |

To understand the algorithm performance varying with problem size, below we increase the number of angular directions to 24, and the total number of unknowns is increased by a factor of 3 to 1,086,090,600. The numerical results with using 1,280, 2,560 5,120 and 10,240 processor cores are summarized in Table III. It is easily found that the number of Newton iterations is kept as a constant for both $MASM_{sub}$ and MASM, and the GMRES iteration of $MASM_{sub}$ is slightly more than that in MASM. We will see that the impact of the increase in GMRES iterations is negligible in the following discussion. The preconditioner setup time of $MASM_{sub}$ is about 10 times smaller than that used in MASM for all core counts. For example, 327 s is taken for MASM at 1,280 cores, while $MASM_{sub}$ takes only 32 s. The preconditioner setup times for $MASM_{sub}$ at 2,560, 5,120 and 10,240 cores are 22 s, 16 s and 19 s, respectively, and these for MASM are 193 s, 185 s, and 168 s. The preconditioner application time for $MASM_{sub}$ is also slightly smaller than MASM for all core counts. The linear solver time and the total compute time of $MASM_{sub}$ are reduced accordingly due to the decrease of the preconditioner setup time and the preconditioner application time, compared with the traditional MASM. The overall algorithm equipped with $MASM_{sub}$ is scalable in terms of the compute time in the sense that the total compute time is halved when we double the number of processor cores. More precisely, the total compute time using $MASM_{sub}$ is 1485 s at 1,280 cores, and it is almost halved to 792 s when we double the number of processor cores to 2,560. It continues being decreased to 443 s and 292 s when the number of processor cores is increased to 5,120 and 10,240, respectively. When the number of processor cores is increased to 10, 240, the parallel efficiency of MASM$_{sub}$ is twice better than that obtained using MASM. When the number of processor cores is 10,240, the total compute time of MASM$_{sub}$ is about 50% of that consumed by MASM. The comparison of the preconditioner setup time and the total compute time between MASM$_{sub}$ and MSM is also shown in Fig. 7, where we, again, observed that the total compute time using MASM$_{sub}$ is much smaller than that consumed by MASM, especially when the number of processor core is large.

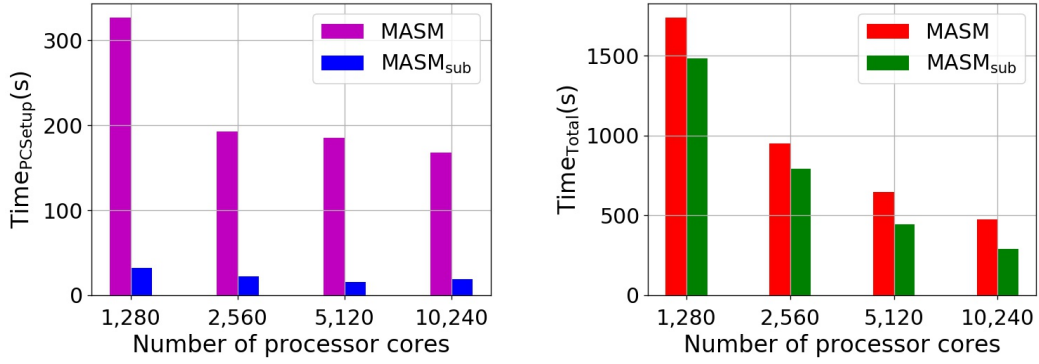

Figure 7. Preconditioner setup time and total compute time using MASM$_{sub}$ and MASM for a problem with 1,086,090,600 unknowns on 1,280, 2,560, 5,120 and 10,240 processor cores. Left: preconditioner setup time, right: preconditioner application time.

Table IV. Performance comparison using 48 angular directions between MASM$_{sub}$ and MASM. The system of eigenvalue equations with 2,172,181,200 unknowns is calculated on 1,280, 2,560, 5,120 and 10,240 processor cores.

| $np$ | PC | Iter$_{Newton}$ | Iter$_{GMRES}$ | Time$_{PCSetup}$ | Time$_{PCApply}$ | Time$_{KSP}$ | Time$_{Total}$ | EFF |
|---|---|---|---|---|---|---|---|---|
| 1,280 | MASM | 7 | 14 | 1022 | 451 | 4079 | 4400 | –% |
| 1,280 | MASM$_{sub}$ | 7 | 15 | 60 | 311 | 3193 | 3513 | 100% |
| 2,560 | MASM | 7 | 14 | 547 | 241 | 2115 | 2280 | 77% |
| 2,560 | MASM$_{sub}$ | 7 | 15 | 37 | 159 | 1612 | 1782 | 99% |
| 5,120 | MASM | 7 | 14 | 398 | 322 | 1476 | 1572 | 56% |
| 5,120 | MASM$_{sub}$ | 7 | 15 | 26 | 86 | 910 | 1006 | 87% |
| 10,240 | MASM | 7 | 14 | 374 | 236 | 1070 | 1128 | 39% |
| 10,240 | MASM$_{sub}$ | 7 | 16 | 19 | 48 | 564 | 621 | 71% |

Next, we further increase the number of angular directions to 48, and the resulting system has 336 unknowns per mesh vertex. The total number of unknowns in the system of eigenvalue equations is 2,172,181,200. The performance comparison between MASM$_{sub}$ and MASM for this configuration is reported in Table IV. Similarly, the numbers of Newton iterations and GMRES iterations are able to stay close to constants when we increase the number of processor cores. The ratio of MASM$_{sub}$ setup time to that in MASM is much higher than that in the previous tests since more unknowns are introduced. That means that a larger time reduction on preconditioner setup is obtained with more unknowns since more preconditioner setup time for MASM is required while the preconditioner setup time of MASM$_{sub}$ does not increase much. More precisely, the full matrix becomes larger with more unknowns, but a submatrix for a given direction of a given energy group stay as the same. Therefore, MASM setup time is significantly increased, while that of MASM$_{sub}$ does not change. At 1,280 cores, 1022 s is taken in MASM setup, while MASM$_{sub}$ uses only 60 s. When we continue increasing the number of processor cores to 2,560, 5,120, and 10,240, the preconditioner setup times for both are decreased, but MASM$_{sub}$ setup is still 20 times

18 FANDE KONG ET AL.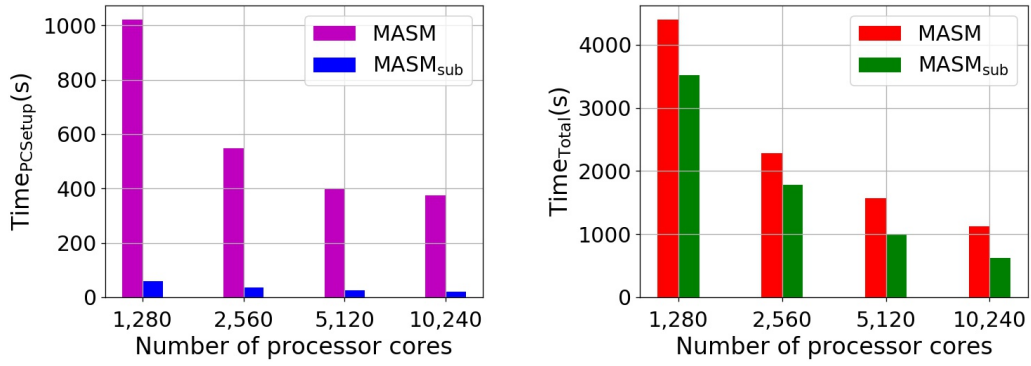

Figure 8. Preconditioner setup time and total compute time for $MASM_{sub}$ and MASM with 48 angular directions on 1,280, 2,560, 5,120, and 10,240 processor cores for a problem with 2,172,181,200 unknowns. Left: preconditioner setup time; right: total compute time.

Table V. Compute times on different components of Newton-Krylov method with 48 angular directions on 1,280, 2,560, 5,120, and 10,240 processor cores for a problem with 2,172,181,200 unknowns.

| $np$ | PC | $Time_{Func}$ | $Time_{Jac}$ | $Time_{LS}$ | $Time_{MF}$ | $Time_{KSP}$ | $Time_{Total}$ | EFF |
|---|---|---|---|---|---|---|---|---|
| 1,280 | MASM | 255 | 65 | 178 | 2626 | 4079 | 4400 | –% |
| 1,280 | $MASM_{sub}$ | 254 | 65 | 177 | 2818 | 3193 | 3513 | 100% |
| 2,560 | MASM | 131 | 34 | 91 | 1329 | 2115 | 2280 | 77% |
| 2,560 | $MASM_{sub}$ | 132 | 39 | 92 | 1415 | 1612 | 1782 | 99% |
| 5,120 | MASM | 73 | 23 | 51 | 792 | 1476 | 1572 | 56% |
| 5,120 | $MASM_{sub}$ | 73 | 23 | 51 | 798 | 910 | 1006 | 87% |
| 10,240 | MASM | 44 | 14 | 31 | 488 | 1070 | 1128 | 39% |
| 10,240 | $MASM_{sub}$ | 44 | 13 | 31 | 497 | 564 | 621 | 71% |

faster than the traditional MASM setup. For the preconditioner application, $MASM_{sub}$ is better than MASM by a factor of 1.45 at 1,280 cores, a factor of 1.5 at 2,560 cores, a factor of 3.7 at 5,120 cores, and a factor of 4.9 at 10,240 cores. The total compute time is improved accordingly for $MASM_{sub}$ thanks to a time reduction in the preconditioner setup. The overall algorithm based on $MASM_{sub}$ is scalable even with up to 10,240 processor cores, where a parallel efficiency of 71% is achieved. The total compute time for $MASM_{sub}$ at 10,240 is half of that in MASM. We conclude that the $MASM_{sub}$ based Newton-Krylov solver is ideally scalable. The preconditioner setup times and the total compute times for $MASM_{sub}$ and MASM are compared in Fig. 8 as well.

To explore more details of Newton-Krylov-$MASM_{sub}$, we below report the compute times on different components of algorithms in Table V. We find that the compute times on Jacobian evaluation, function evaluation, line search, and matrix-free operation are the same for both $MASM_{sub}$ and MASM, and are also scalable for all core counts. The matrix-free operation accounts for 50% of the total compute time since each GMRES iteration based on the matrix-free method involves one function evaluation at the current solution.

Again, we continue increasing the number of angular directions to 80, and the resulting system has 560 variables per mesh vertex. We report numerical results for this problem with 3,620,302,000 unknowns in Table VI. There is no data on 1,280 processor cores since the problem is too large and the available memory is not enough. From Table VI, the number of Newton iterations is 7 for $MASM_{sub}$ and MASM for all core counts. The number of GMRES iterations for $MASM_{sub}$ is similar to that of MASM. The preconditioner setup of $MASM_{sub}$ is 20 times more efficient than the traditional MASM setup for all core counts. For example, MASM uses 1172 s at 2,560 cores, while that of $MASM_{sub}$ is 62 s. The preconditioner application time of $MASM_{sub}$ is also much smaller than that used in MASM. Due to the time reduction on both the precondtioner setup and the preconditioner application, the total compute time is significantly decreased, compared with MASM. When the number of procesor cores is 10,240, the total compute time of $MASM_{sub}$ is only half of that taken by MASM, and the corresponding parallel efficiency is twice better. In all, the new proposed algorithm is able to maintain an ideal speedup with up to 10,240 processor cores. The performance comparison on the preconditioner setup time and the total compute time between $MASM_{sub}$ and MASM is also added in Fig. 9, where the same behavior is observed.

Table VI. Parallel performance with 80 angular directions using 2,560, 5,120, and 10,240 processor cores. A system of eigenvalue equations with 3,620,302,000 unknowns is computed using Newton-Krylov together with $MASM_{sub}$ and MASM.

| $np$ | PC | $Iter_{Newton}$ | $Iter_{GMRES}$ | $Time_{PCSetup}$ | $Time_{PCApply}$ | $Time_{KSP}$ | $Time_{Total}$ | EFF |
|---|---|---|---|---|---|---|---|---|
| 2,560 | MASM | 7 | 14 | 1172 | 448 | 4398 | 4723 | –% |
| 2,560 | $MASM_{sub}$ | 7 | 15 | 62 | 267 | 3592 | 3930 | 100% |
| 5,120 | MASM | 7 | 14 | 760 | 498 | 2603 | 2788 | 70% |
| 5,120 | $MASM_{sub}$ | 7 | 15 | 39 | 140 | 1663 | 1846 | 106% |
| 10,240 | MASM | 7 | 14 | 586 | 411 | 1893 | 2006 | 49% |
| 10,240 | $MASM_{sub}$ | 7 | 16 | 30 | 79 | 1064 | 1176 | 84% |

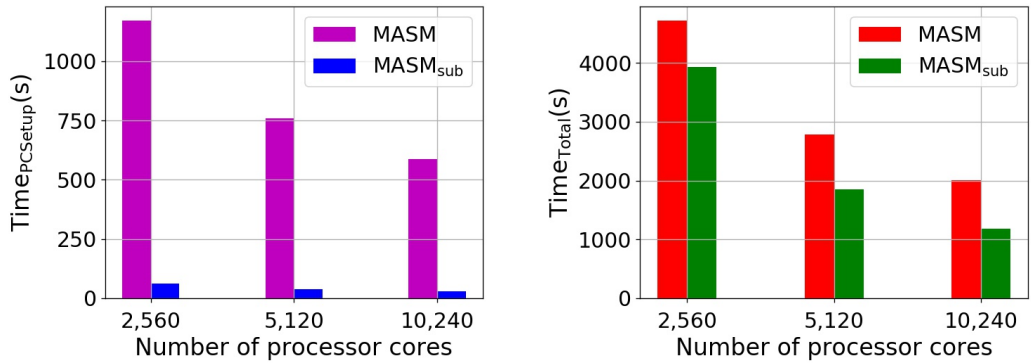

Figure 9. Preconditioner setup time and total compute time using 80 angular directions on 2,560, 5,120, and 10,240 processor cores for a problem with 3,620,302,000 unknowns. Left: preconditioner setup time; right: total compute time.

Now, let us summarize our discussions. The first thing we concern on is the variaion of the preconditioner setup time with problem size, and it is summarized in Fig. 10 using 2,560 and 10,240 processor cores. We observed, in Fig 10, that when we increase the number of unknowns per mesh vertex, the $MASM_{sub}$ setup time does not increase much, while that of MASM is significantly increased. As we stated earlier, it is because the full matrix MASM uses is proportionately increased when the number of unknowns is increased, while the submatrix size for $MASM_{sub}$ does not change. The overall parallel efficiencies and speedups are summarized in Fig. 11 for all problem sizes. Finally, we conclude that the overall algorithm based on $MASM_{sub}$ is capable of maintaining an ideal scalability as long as the problem size is not too small. The good scalabilities using 24, 48, 80 angular directions, respectively, are achieved with up to 10,240 processor cores.

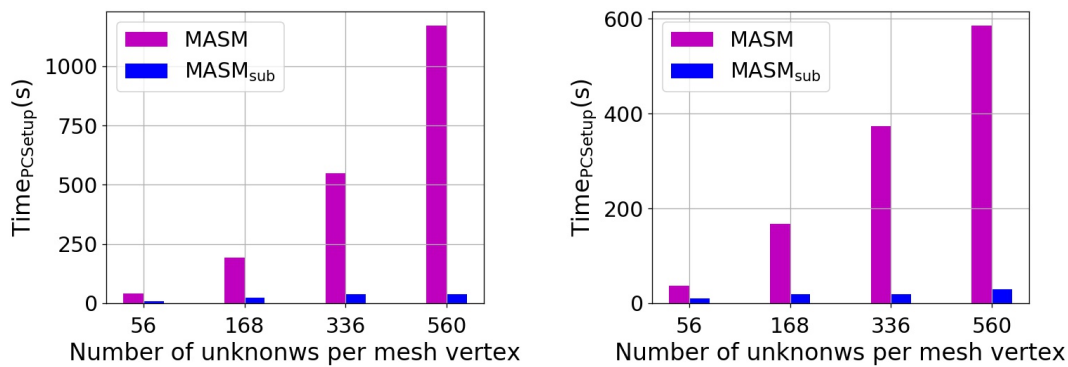

Figure 10. Preconditioner setup time varying with problem size. Left: preconditioner setup time using 2,560 processor cores for different problem sizes; right: preconditioner setup time using 10,240 processor cores for different problem sizes.

## 5. CONCLUSIONS

We have studied a scalable parallel multilevel domain decomposition preconditioner for solving the algebraic eigenvalue system arising from the discretization of the multigroup neutron transport equations. The multigroup neutron transport equations is discretized using the first-order finite element method in space and using the discrete ordinates method in angle. The fully coupled Newton-based eigenvalue solver was employed to compute the largest eigenvalue of the neutron transport problem. During each Newton iteration, the Jacobian system was calculated using GMRES together with the multilevel methods. It is often expensive to construct necessary coarse spaces for the multilevel method, and thus we introduced an inexpensive approach, the subspace-based coarsening algorithm, to construct coarse spaces for the proposed parallel algorithm framework. Compared with the traditional multilevel methods, the proposed approach is much faster on the construction of coarse spaces especially when the number of unknowns per mesh vertex is large. Due to the time reduction on the preconditoner setup, the overall algorithm is shown to be scalable on a large number of processor cores. We have numerically demonstrated that the proposed algorithm is able to maintain an ideal scalability with more than 10,000 processor cores for the 3D C5G7 benchmark on unstructured meshes with billions of unknowns.

## ACKNOWLEDGMENTS

This manuscript has been authored by Battelle Energy Alliance, LLC under Contract No. DE-AC07-05ID14517 with the U.S. Department of Energy. The United States Government retains and the publisher, by accepting the article for publication, acknowledges that the United States Government retains a nonexclusive, paid-up, irrevocable, and worldwide license to publish or reproduce the published form of this manuscript, or allow others to do so, for United States Government purposes.

This research made use of the resources of the High-Performance Computing Center at Idaho National Laboratory, which is supported by the Office of Nuclear Energy of the U.S. Department of Energy and the Nuclear Science User Facilities under Contract No. DE-AC07-05ID14517.

REFERENCES

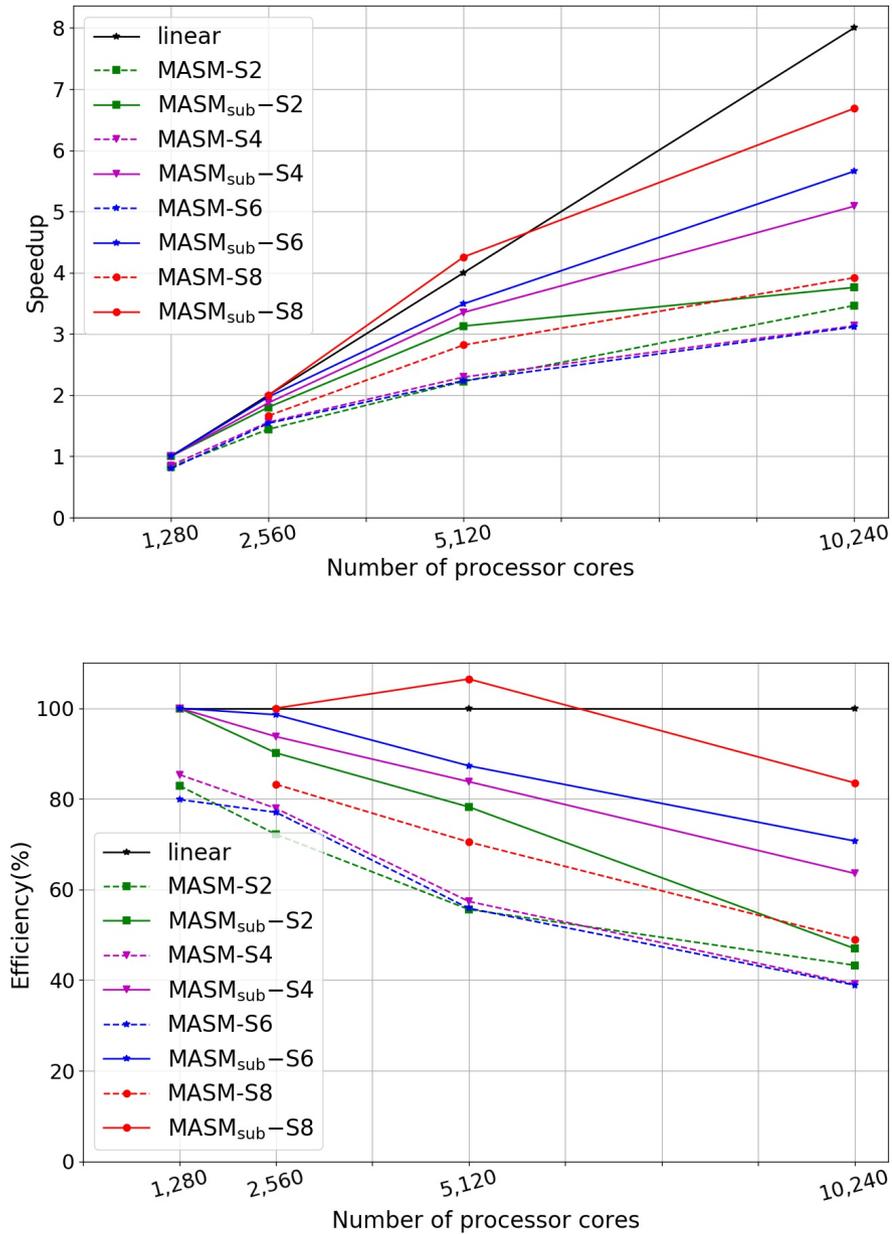

Figure 11. Speedups and parallel efficiencies for all problem sizes using 1,280, 2,560, 5,120, and 10,240 processor cores. Top: speedups, bottom: parallel efficiencies. "Sx" represents the number of unknowns per mesh vertex. "S2": 56; "S4": 168; "S6": 336; "S8": 560.